\documentclass{amsart}
\newcommand{\I}{\mathcal I}
\newcommand{\IM}{{\mathcal I}^{\max}}
\newcommand{\N}{\mathbb N}
\newcommand{\Z}{\mathbb Z}
\newcommand{\del}{\!\setminus\!}
\newtheorem{theorem}{Theorem}[section]
\newtheorem{lemma}[theorem]{Lemma}

\begin{document}
\title{Relative rank axioms for infinite matroids}
\author{R.A. Pendavingh}
\begin{abstract} In a recent paper, Bruhn, Diestel, Kriesell and Wollan present four systems of axioms for infinite matroids, in terms of independent sets, bases, closure and circuits. No system of rank axioms is given. We give an easy example showing that rank function of an infinite matroid may not suffice to characterize it. We present a system of axioms in terms of {\em relative rank}.\end{abstract}
\date{\today}
\maketitle
\section{Introduction}
In \cite{BDKW}, an infinite matroid is defined as a pair $M=(E,\I)$ where $E$ is any set, and $\I$ is a set of subsets of $E$, satisfying the following {\em independence axioms}:
\begin{enumerate}
\item[(I1)] $\emptyset\in \I$.
\item[(I2)] $\I$ is closed under taking subsets.
\item[(I3)] if $I\in \I\del \IM$ and $I'\in \IM$, then $I+x\in \I$ for some $x\in I'\del I$.
\item[(IM)] if $I\subseteq X\subseteq E$ and $I\in \I$, then the set $\{I'\in \I\mid I\subseteq I'\subseteq X\}$
has an inclusionwise maximal element.
\end{enumerate}
Here, $\IM$ denotes the set of (inclusionwise) maximal elements in $\I$. 

In \cite{BDKW}, three further axiom systems are described that equivalently define infinite matroids in terms of bases, circuits, and a closure operator respectively. There is no axiomatization in terms of the rank function, and it is not hard to see why such a description infinite matroids cannot be given. The rank of a set $X$ in an infinite matroid is inevitably 
$$r(X):=\max\{ |I| \mid I\in \I, I\subseteq X\}.$$ 
But then the infinite matroids $M=(\Z, 2^{\Z})$ and $M'=(\Z, 2^{\Z}\del\{\Z\})$ have identical rank functions. So the rank function does not even suffice to characterise an infinite matroid, in general.

In a finite matroid, we could define the {\em relative rank} of an extension $A\supseteq B$ as $r(A|B):=r(A)-r(B)$. From this definition and the rank axioms, one easily derives the following.
\begin{enumerate}
\item[(R1)] $0\leq r(A|B)\leq |A\del B|$ for all $B\subseteq A\subseteq E$.
\item[(R2)] $r(A|A\cap B)\geq r(A\cup B|B)$ for all $A, B\subseteq E$.
\item[(R3)] $r(A|C)=r(A|B)+r(B|C)$ for all $C\subseteq B\subseteq A\subseteq E$.
\end{enumerate}
It is not difficult to see that in turn, the rank axioms follow for $r(X):=r(X|\emptyset)$ from these three properties. 

In this note we show that the relative rank can be extended to infinite matroids in such way that an axiomatization of infinite matroids in terms of relative rank is possible. 
 
\section{Relative rank in Infinite matroids}
In the infinite matroid $M=(E,\I)$, we may define the {\em relative rank} 
of a pair $A, B\subseteq E$ so that $A\supseteq B$ as 
$$r_M(A|B):=\max\{ |I\del J| \mid J\subseteq I\subseteq A,  ~I\in\I, ~J\text{ max. independent in } B\}.$$
In the following lemmas we assume that $M=(E, \I)$ is a fixed infinite matroid. The next lemma shows that $r_M$ is well-defined.
\begin{lemma} \label{witness} For any $B\subseteq A\subseteq E$ there exist $I, J\in \I$ so that $I\supseteq J$, $I$ is maximally independent within $A$ and $J$ is maximally independent within $B$. For any such $I, J$, we have $r_M(A|B)=|I\del J|$.\end{lemma}
\proof By Theorem 3.4 of \cite{BDKW}, the restriction $M|A:=(A, \I\cap 2^A)$ is a matroid. So we may assume $E=A$. By (IM), there is a maximal independent set $J$ in $B$ and a maximally independent set $J'$ in $A$. Again by (IM), there is a maximal set $I$ in $\{I\in\I\mid J\subseteq I\subseteq J\cup J'\}$. By (I3), $I$ is maximally independent in $A$, as required. To see that $r_M(A|B)=|I\del J|$, let $I',J'$ attain the maximum in the definition of rank. By Lemma 3.5 of \cite{BDKW}, both $I\del J$ and $I'\del J'$ are max. independent in $M/B$. Then by Lemma 3.7 of \cite{BDKW}, we have $|I\del J|=|I'\del J'|$.
\endproof 
The next lemma allows us to zoom in on a minor of $M$, which will be convenient in what follows.
\begin{lemma} \label{zoom}Let $M=(E,\I)$ be an infinite matroid, let $X,Y\subseteq E$ be disjoint sets and let $N=(M/X)|Y$.
Then for any $A, B\subseteq E$ so that 
$X\subseteq B\subseteq A\subseteq Y$, we have 
$r_M(A|B)=r_N(A\del X|B\del X)$.
\end{lemma}
\proof In the special case that  $X=B$ and $A=Y$, this lemma follows from the previous one. Otherwise, we have 
$$r_M(A|B)=r_{(M/B)|A}(A\del B, \emptyset)=r_N(A\del X|B\del X)$$
by two applications of the special case.
\endproof
We next show five properties which together will turn out to characterise relative rank functions. The first follows directly from the definition of $r_M$.
\begin{lemma} $0\leq r_M(A|B)\leq |A\del B|$.\end{lemma}
\begin{lemma} $r_M(A|A\cap B)\geq r_M(A\cup B|B)$ for all $A, B\subseteq E$.\end{lemma}
\proof By Lemma \ref{zoom}, we may assume that $A\cap B=\emptyset$. By Lemma \ref{witness}, there exist $I,J$ so that $J\subseteq I$, $I$ is maximally independent in $A\cup B$ and $J$ is maximally independent in $B$. Then $r_M(A|A\cap B)=r_M(A\cup B,\emptyset)\geq |I\del J|=r_M(A\cup B|B)$.\endproof
\begin{lemma}$r_M(A|C)=r_M(A|B)+r_M(B|C)$ for all $C\subseteq B\subseteq A\subseteq E$.\end{lemma}
\proof By Lemma \ref{zoom}, we may assume $C=\emptyset$. By Lemma \ref{witness}, there exist $I,J$ so that $J\subseteq I$, $I$ is maximally independent in $A$ and $J$ is maximally independent in $B$. Then $r_M(A|C)=|I|=|I\del J|+|J|=r_M(A|B)+r_M(B|C)$.\endproof
\begin{lemma} if $A=\bigcup_{\gamma\in \Gamma} A_\gamma$, and $r_M(A_\gamma|B)=0$ for all $\gamma\in\Gamma$, then $r_M(A|B)=0$.\end{lemma}
\proof By Lemma \ref{zoom}, we may assume $B=\emptyset$. If $r_M(A_\gamma|B)=0$ for all $\gamma\in\Gamma$, then there is no independent singleton in $A$. Hence $r_M(A|B)=0$.\endproof
\begin{lemma}\label{RM} for all $A, B\subseteq E$ so that $A\supseteq B$, there exist an $I\in\I$ so that $r_M(A| I)=0$ and $r_M(B| B\cap I)=0$.\end{lemma}
\proof Let $I, J$ be as in Lemma \ref{witness}. If $I'$ is an independent set in $M/I$, then $I\cup I'$ is independent in $M$ by Lemma 3.5 of \cite{BDKW}, so by maximality of $I$ and Lemma \ref{zoom} we have $r_M(A|I)=0$. Similarly, it follows that $r_M(B|J)=0$ and $J\supseteq B\cap I$.\endproof

Knowing $r_M$ suffices to recover $\I$.
\begin{lemma}\label{recover}$I\in \I$ if and only if $r_M(I, I-x)>0$ for all $x\in I$.\end{lemma}
\proof Necessity is straightforward. To see sufficiency, suppose $I$ is not independent, and contains a maximal independent set $J\subseteq I$. Then $r_M(I|J)=0$, hence $r_M(I|I-x)=0$ for any $x\in I\del J$.\endproof

The rank functions of a matroid and its dual have an easy relationship, which in fact characterizes the dual.
\begin{lemma} Let $M$ and $M'$ be infinite matroids with common ground set $E$. Then $M'=M^*$ if and only if
$$r_M(A|B)+r_{M'}(E\del B|E\del A)=|A\del B|$$ for all $B\subseteq A\subseteq E$.\end{lemma}
\proof Necessity: suppose $M'=M^*$, and consider $B\subseteq A$. We may assume $B=\emptyset$, $A=E$. If $I$ is a basis of $M$, then  $E\del I$ is a basis of $M^*$, hence $r_M(A|B)+r_M^*(E\del B|E\del A)=|I|+ |E\del I|=|E|=|A\del B|$.

Sufficiency: suppose $M'\neq M^*$. Then by Lemma \ref{recover} we have 
$r_{M^*}(A|B)\neq r_{M'}(A|B)$ for some $A,B$. By Lemma \ref{zoom}, we may assume $B=\emptyset$ and $A=E$. In case $r_{M^*}(E|\emptyset)< r_{M'}(E|\emptyset)$,  a maximal independent set of $M^*$ has $r_{M^*}(E|I)=0< r_{M'}(E|I)$. Then 
$r_{M^*}(I+x|I)=0< r_{M'}(I+x|I)$ for some $x$, and again by Lemma \ref{zoom} we may assume $E=\{x\}$. Then $r_M(E|\emptyset)+r_{M'}(E|\emptyset)=0\neq 1=|E|$, as required.
The case that  $r_{M^*}(E|\emptyset)> r_{M'}(E|\emptyset)$ is similar.
\endproof

\section{Relative rank axioms for infinite matroids}
We consider partial functions
$$r:2^E\times 2^E\rightarrow \N\cup\{\infty\}$$
so that $r(A|B)$ is defined if $B\subseteq A\subseteq E$. 
For such an $r$, we define 
$$\I_r:=\{I\subseteq E\mid r(I| I-x)>0\text{ for all }x\in I\}$$
and we will say that $I$ is {\em $r$-independent} if $I\in \I_r$.
\begin{theorem} Let $E$ be a set, and let $r:2^E\times 2^E\rightarrow \N\cup\{\infty\}$ be a partial function such that 
\begin{enumerate}
\item[(R1)] $0\leq r(A|B)\leq |A\del B|$ for all $B\subseteq A\subseteq E$
\item[(R2)] $r(A|A\cap B)\geq r(A\cup B|B)$ for all $A, B\subseteq E$
\item[(R3)] $r(A|C)=r(A|B)+r(B|C)$ for all $C\subseteq B\subseteq A\subseteq E$
\item[(R4)] if $A=\bigcup_{\gamma\in \Gamma} A_\gamma$, and $r(A_\gamma|B)=0$ for all $\gamma\in\Gamma$, then $r(A|B)=0$
\item[(R5)] for all $A, B\subseteq E$ so that $A\supseteq B$, there exist an $I\in\I_r$ so that $r(A| I)=0$ and $r(B| B\cap I)=0$ 
\end{enumerate}
Then $M=(E, \I_r)$ is an infinite matroid, and $r=r_M$.
\end{theorem}
\proof The proof takes the form of a number of lemmas. Trivially, we have:
\begin{lemma} \label{I1} For any $r$,  we have 
\begin{enumerate}\item[(I1)] $\emptyset\in\I_r$.
\end{enumerate}
\end{lemma} 
The submodularity of $r$ suffices to show that $\I_r$ is closed under taking subsets. 
\begin{lemma} \label{I2}If $r$ satisfies (R2), then  
\begin{enumerate}\item[(I2)] $\I_r$ is closed under taking subsets.\end{enumerate}
\end{lemma} 
\proof Let $I\in \I_r$, and let $J\subseteq I$. If $J\not \in \I_r$, then $r(J| J-x)=0$ for some $x\in J$. Taking $A=J, B=I-x$ in (R2), we find $r(J|J-x)\geq r(I|I-x)$. Hence $r(I|I-x)=0$ and $I\not\in\I_r$, a contradiction.
\endproof
\begin{lemma} \label{plus}Let $r$ satisfy (R1), (R3). If $I\in \I_r$, then 
$$ \label{Icrit}I+x\in\I_r \Longleftrightarrow r(I+x|I)>0.$$
\end{lemma}
\proof Necessity is immediate from the definition of $\I_r$. To see sufficiency, note that if $I+x\not\in \I_r$, there is a $y$ so that $r(I+x|I+x-y)=0$. Then by (R1) and (R3), we have
$$r(I+x|I-y)=r(I+x|I+x-y)+r(I+x-y|I-y)\leq 1.$$
Also by (R3), we have 
$$r(I+x|I-y)=r(I+x|I)+r(I|I-y).$$
Since $I$ is assumed independent, we have $r(I|I-y)>0$. Hence $r(I+x|I)=0$, as required. \endproof

\begin{lemma} \label{span}Let $r$ satisfy (R1), (R3), and (R4). If $I\in \I_r$ and $I\subseteq F\subseteq E$, then 
$$ I\text{ is maximally $r$-independent in }F \Longleftrightarrow r(F|I)=0.$$
\end{lemma}
\proof Necessity: by Lemma \ref{plus}, we have $r(I+x|I)=0$ for all $x\in F\del I$.  By applying (R4) to $\Gamma=F\del I$, $A_x=I+x$, $B=I$, we find that $r(F|I)=0$. 

Sufficiency: if $I$ is not maximally $r$-independent in $F$, then $r(I+x|I)>0$ for some $x\in F\del I$, and then $r(F|I)=r(F|I+x)+r(I+x|I)>0$ by (R3). \endproof

\begin{lemma} Let $r$ satisfy (R1), (R3), (R4). Then 
\begin{enumerate}\item[(I3)] if $I\in \I_r\del \IM_r$ and $I'\in \IM_r$, then $I+x\in \I_r$ for some $x\in I'\del I$. \end{enumerate} 
\end{lemma}
\proof Suppose $I\in \I_r\del \IM_r$ and $I'\in \IM_r$, and that $I+x\not\in \I_r$ for all $x\in I'\del I$. Taking $F=I\cup I'$ in Lemma \ref{span}, we find that $r(I\cup I'| I)=0$. Since $I'$ is maximally independent, we have $r(E, I')=0$, so that $r(E|I\cup I')=0$ by (R3). Again by (R3), we have $r(E|I)=r(E|I\cup I')+r(I\cup I'|I)=0$, so that $I$ is maximally independent, a contradiction.
\endproof

\begin{lemma} Let  $r$ satisfy (R1), (R3), (R4), and (R5). Then 
\begin{enumerate}\item[(IM)] if $I\subseteq X\subseteq E$ and $I\in \I_r$, then the set $\{I'\in \I_r\mid I\subseteq I'\subseteq X\}$
has an inclusionwise maximal element. \end{enumerate} 
\end{lemma}
\proof Let $I\in\I_r$ and $X\supseteq I$. Applying (R5) with $A=X$ and $B=I$, we find a $J\in \I_r$ so that $r(X|J)=0$ and $r(I| I\cap J)=0$. As $I\in\I_r$, we have  $I\subseteq J$, and by Lemma \ref{span}, $J$ is a maximal element of  $\{I'\in\I_r\mid I\subseteq I'\subseteq Y\}$. \endproof
To finish the proof, let us assume that $r$ satisfies all the axioms. Then the above lemmas establish that $M:=(E, \I_r)$ is an infinite matroid. It remains to show that $r$ is the relative rank function of $M$. If not, then $r_M(A|B)\neq r(A|B)$ for some $A, B$, and without loss of generality $A=E$ and $B=\emptyset$. Consider a maximal $r$-independent set $I$. Then $r_M(A|B)=|I|$ by definition, and $r(A|I)=0$ by Lemma \ref{span}, so that $r(A|B)=r(I|\emptyset)$. But for an $r$-independent set $I$, we have $r(I|\emptyset)=|I|$ by induction, as $r(I|I-x)=1$ and $I-x$ is again $r$-independent. This completes the proof of the Theorem.\endproof
As is evident from the proof, axioms (R1)---(R4) imply (I1), (I2), (I3), so that replacing (R5) with 
\begin{enumerate}
\item[(RM)] $\I_r$ satisfies (IM)
\end{enumerate}
would give an equivalent system of axioms. Such a system of axioms would be more in line with axioms as given in Section 1 of \cite{BDKW}. 

Finally, we note that if $r(E|\emptyset)$ is finite, then (R5) follows from the other axioms, and if $E$ is a finite set, then (R4) is redundant as well. 


\begin{thebibliography}{1}
\bibitem{BDKW} H. Bruhn, R. Diestel, M. Kriesell \& P. Wollan, {\em Axioms for infinite matroids}, preprint arXiv:1003.3919 (2010) . 
\end{thebibliography}
\end{document}